\documentclass[12pt]{article}
\usepackage{amsmath, amsfonts, amssymb, amsthm} \parskip1pt
\usepackage{mathrsfs}
\usepackage{amsfonts,graphics,epsfig}
\usepackage{amsmath, amsfonts, amssymb, amsthm, amscd, setspace, graphicx}
\newcommand{\PSL}{\operatorname{PSL}}
\newcommand{\SL}{\operatorname{SL}}
\newtheorem{theorem}{Theorem}[section]
\newtheorem{proposition}[theorem]{Proposition}
\newtheorem{definition}[theorem]{Definition}
\newtheorem{remark}[theorem]{Remark}

\newtheorem{lemma}[theorem]{Lemma}

\def\R{\mathbb R} \def\Z{\mathbb Z}
 \def\H{\mathbb H}

\def\N{\mathbb N}
 \def\d{{\rm d}}

\def\N{{\bf N}}

\date{}

\def\disp{\displaystyle}

\def\build#1_#2^#3{\mathrel{\mathop{\kern 0pt#1}\limits_{#2}^{#3}}}

\def\smallsquare{\vbox{\hrule\hbox{\vrule height 1 ex\kern 1
ex\vrule}\hrule}}

\def\N{{\mathbb N}}

\begin{document}
\title{On certain orbits of geodesic flow and $(a,b)$-continued fractions}

\author{ Manoj Choudhuri,\\ 
 {\it Institute of Infrastructure Technology Research and Management},\\ 
{\it Near Khokhara Circle, Maninagar (East),}\\ {\it Ahmedabad-380026, Gujarat, India}\\
email: manojchoudhuri@iitram.ac.in}

\maketitle
\begin{abstract}
In this article, we characterize two kinds of exceptional orbits of the geodesic flow associated with the Modular surface in terms of a two-parameter family of continued
fraction expansion of endpoints of the lifts to the hyperbolic plane of the corresponding geodesics. As a consequence, we obtain an extension of Dani correspondence between 
homogeneous dynamics and Diophantine approximation.
 \vspace*{0.2cm}
 
 

\end{abstract}
\section{Introduction}
Let $\H := \{z=x+iy:y>0\}$ be the upper half plane endowed with the hyperbolic metric $ds^2=\frac{\disp{dx^2+dy^2}}{\disp{y^2}}$. The group $\PSL(2,\R)=\SL(2,\R)/\{\pm I\}$ acting by fractional linear transformations (see \cite{K}), is the orientation preserving isometry group of $\H$. The discrete group $\PSL(2,\Z)=\SL(2,\Z)/\{\pm I\}$ acting properly discontinuously on $\H$ gives rise to the modular surface $M=\H/\PSL(2,\Z)$, which is topologically a sphere with two singularities and one cusp. Let $T^1\H$ be the unit tangent bundle of the hyperbolic plane which is the collection $(z,\zeta)$ with $z$ in $\H$ and $\zeta$ being a tangent vector of norm one in $T_z\H$; $\PSL(2,\Z)$ acts on $T^1\H$ as well and the quotient space $T^1\H/\PSL(2,\Z)$ can be identified with the unit tangent bundle of $M$, which we denote by $T^1M$. We denote $(z,\zeta)\in T^1M$ by $v$. Any $v\in T^1M$ determines a unique geodesic in $M$. If we consider the geodesic along with its tangent vector at each point, then it is the orbit of $v$ under the geodesic flow. This orbit is denoted by $\{g_t v\}$, where $g_t$ denotes the geodesic flow on $T^1M$. It is well known that the geodesic flow on $T^1M$ is ergodic with respect to the Liouville measure (see \cite{EW} for more details and original references). 
This means in particular that with respect to the Liouville measure, the orbits of almost all $v$ in $T^1M$ are equidistributed, i. e., if $\mu$ denotes the (normalized) Liouville measure 
on $T^1M$, then for almost all $v\in T^1M$, $$\frac{1}{T}\int_0^T \chi_A(g_tv)\underset{T\longrightarrow\infty}{\xrightarrow{\hspace*{1.2cm}}} \mu(A),$$  for any measurable $A\subset T^1M$, where $\chi_A$ denotes the characteristic function of the set $A$. Apart from these generic orbits there are many interesting orbits of geodesic flow associated with the modular surface. By Dani correspondence (see \cite{D}, \cite{D1} for details), we know that badly approximable numbers (see \cite{HW} for definition) correspond to bounded orbits and rational numbers correspond to divergent orbits. In this article, we are going to characterize two kinds of orbits in terms of their asymptotic rate of time spent in cusp neighbourhoods. One of these kinds of orbits contains the bounded orbits and the other one contains the divergent orbits. These characterizations are done in terms of the continued fraction expansions of certain real numbers associated with those orbits. In order to do so, we use arithmetic coding of geodesics on the modular surface which was originated in the 1924 paper of E. Artin (\cite{A}), who proved the existence of a dense geodesic using classical continued fraction. A more precise description of this machinery using classical continued fraction can be found in \cite{S}. We do not restrict ourselves to only classical continued fraction, rather we consider a two-parameter family of continued fractions and obtain the characterizations in terms of that family, from which the characterizations in terms of the classical continued fraction follow. 

The two-parameter family of continued fractions we are going to consider in this article is known as $(a,b)$-continued fractions with $a$, $b\in\R$ satisfying a technical condition (see next section for more details about $(a,b)$-continued fractions). Using these continued fraction expansions of real numbers, in \cite{KU5}, S. Katok and I. Ugarcovicci describe a coding of geodesics on the modular surface, which enables one to give a symbolic description of the geodesic flow associated with the modular surface. We use this coding of geodesics, which also relies on another paper (\cite{KU4}) by the same authors, as a base of the arguments, used to prove the results of the present article. In this article, we consider $(a,b)$-continued fractions for $(a,b)$ in a particular subset $\mathcal{P}$ of $\R^2$, where $\mathcal{P}$ is given as follows:  
\begin{center}
 $\mathcal{P}=\{(a,b)\in\R^2|-1\leq a< 0 < b\leq1,b-a\geq 1\}$.
\end{center}
Note that this excludes the possibilities $a<-1$ and $b>1$, though in the work of Katok and Ugarcovicci (\cite{KU4}, \cite{KU5}) those possibilities were also considered with $-ab\leq1$. Also let $\mathcal{E}$ be the exceptional set discussed in \cite{KU4}, the elements of which do not satisfy the finiteness 
(see next section for the definition) condition. Let $$\mathfrak{S}:=\{(a,b)\in\mathcal{P}\backslash\mathcal{E}:a,b \hspace*{0.1cm}\text{have strong cycle property}\},$$ (see next section for definition of cycle properties) and  
$$\mathcal{S}=\mathfrak{S}\cup \left\{(-1,-1),\left(-\frac{1}{2},\frac{1}{2}\right)\right\}.$$
Now let $\H_{\displaystyle d}:=\{x+iy\in \H:y>d\}$ and $\overline{\H}_{\displaystyle d}\subset T^1\H$ be given by $$\overline{\H}_{\displaystyle d}:=\bigcup\limits_{\mathcal{I}(z)>d}T^1_z\H,$$ where $\mathcal{I}(z)$ denotes the imaginary part of the complex number $z$, and $T^1_z\H$ denotes the set of unit tangent vectors in $T_z\H$. Let $\pi$ denote both the projections from $\H$ to $M$ and from $T^1\H$ to $T^1M$. Also let 
$M_{\displaystyle d}:=\pi(\H_{\displaystyle d})$ and $\overline{M}_{\displaystyle d}
:=\pi(\overline{\H}_{\displaystyle d})$. Note that $\overline{M}_{\displaystyle d}$ is a typical neighbourhood of the cusp in $T^1M$. We say that an orbit $\{g_tv\}_{t\geq0}$ visits the cusp with {\it frequency} $0$, if there exists some $d>1$ such that $\frac{\disp{1}}{\disp{T}}\disp{\int_0^T\chi_
{\overline{M}_{\displaystyle d}}(g_tv)dt}\rightarrow 0$ as $T\rightarrow\infty$. On the other hand, an orbit $\{g_tv\}_{t\geq0}$ is said to visit the cusp with {\it frequency} $1$, if for all $d>1$, $\frac{\disp{1}}{\disp{T}}\disp{\int_0^T\chi_{\overline{M}_{\displaystyle d}}(g_tv)dt}\rightarrow 1$ as $T\rightarrow\infty$. Now let $x\in\R$ and 
$x:=[a_0,a_1,...]_{a,b}$ be its $(a,b)$-continued fraction expansion. Given $\xi>1$ and $j\geq0$, 
we define the modified partial quotients of the $(a,b)$-continued fraction expansion of $x$ as follows:
\begin{align*}
a_j^{(\xi)} &= a_j, \text{if} \hspace{0.2cm} |a_j|>\xi, \\
a_j^{(\xi)} &= 1, \text{if} \hspace{0.2cm} |a_j|\leq \xi.
\end{align*}
\begin{theorem}\label{Cusp 1.1}
For a given $v\in T^1M$, let $\gamma_v$ be the corresponding geodesic in $M$, and $\tilde{\gamma}_v$ be one of its lifts to the hyperbolic plane. Let $x$ be the attracting end point of $\tilde{\gamma_v}$. For $(a,b)\in \mathcal{S}$, let the $(a,b)$-continued fraction expansion of $x$ be given
by $x=[a_0,a_1,a_2,...]_{a,b}$. Also for $\xi>1$, let $\{a_j^{(\xi)}\}_{j\geq0}$ be the modified sequence of partial quotients as defined above, and
$$A_N^{(\xi)}=\frac{\displaystyle1}{\displaystyle N}\displaystyle{\sum\limits_{\disp{j=0}} ^{\disp{N-1}}}\log|a_j^{(\xi)}|, \hspace*{0.3cm}   A_N=\frac{\displaystyle1}{\displaystyle N}\displaystyle{\sum\limits_{\disp{j=0}}^{\disp{N-1}}}\log|a_j|.$$
Then,\\\\ $(i)$ the forward orbit $\{g_t v\}_{t\geq0}$  visits the cusp with frequency $0$ if and only 
if\\ $A_N^{(\xi)}\rightarrow 0$ as $N\rightarrow\infty$ for some $\xi>1$, and\\\\
$(ii)$ $\{g_t v\}_{t\geq0}$ visits the cusp with frequency $1$ if and
only if $A_N\rightarrow\infty$ as $N\rightarrow\infty$. 
\end{theorem}
\begin{remark}
The restriction of the parameters $(a,b)$ to the set $\mathcal{S}$ ensures that any geodesic in $\H$ is $\PSL(2,\Z)$-equivalent to an $(a,b)$-reduced geodesic (a notion to be made clear in the next section). This is essential for this article as we are not looking at a generic set of orbits of the geodesic flow. It follows from Theorem $7.1$ of \cite{KU4} that if $a$ and $b$ have the strong cycle property, then every geodesic in $\H$ is $\PSL(2,\Z)$-equivalent to an $(a,b)$-reduced geodesic, whereas in \cite{KU}, the same statement was shown to be true for $(a,b)=(-1,1),(-\frac{1}{2},\frac{1}{2})$. It is easy to see that $(-1,1)$-continued fraction expansion of a real number $x$ is nothing but the classical continued fraction expansion of $x$ with alternating signs (see \cite{KU} for details). A similar relation holds between $(-\frac{1}{2},\frac{1}{2})$-continued fraction expansion and the 
nearest integer continued fraction expansion of a real number. Then 
the characterizations in Theorem \ref{Cusp 1.1} in terms of the classical and the nearest integer continued fractions follow from the characterizations in terms of $(-1,1)$ and $(-\frac{1}{2},\frac{1}{2})$-continued fractions respectively.
\end{remark}

If we consider the algebraic description of the geodesic flow, then Theorem \ref{Cusp 1.1} may be thought of as an extension of Dani correspondence between homogeneous dynamics and Diophantine approximation. We know that $\PSL(2,\R)=\SL(2,\R)/\{\pm I\}$ can be identified with $T^1\H$ (see \cite{K} for details), where the 
identification is given by 
\begin{center}
$g\longmapsto g(i,\hat{i})$
\end{center}
for $g\in\PSL(2,\R)$, here $\hat{i}$ denotes the unit tangent vector based at the point $i$ and pointing upwards.
Similarly $\PSL(2,\Z)\backslash\PSL(2,\R)\simeq\SL(2,\Z)\backslash\SL(2,\R)$ can be identified with $T^1M$. The right action of the one parameter subgroup\\ $\left\{ a_t:=
\begin{pmatrix}
e^{-\frac{t}{2}} & 0\\
0& e^{\frac{t}{2}}
\end{pmatrix}
\right\}$ on $\SL(2,\Z)\backslash\SL(2,\R)$, which is given by the following:
\begin{center}
$(\Gamma g)a_t\mapsto(\Gamma)ga_t$, for $g\in\SL(2,\R)$,
\end{center}
where we denote $\SL(2,\Z)$ by $\Gamma$, corresponds to the geodesic flow on $T^1M$. Given a real number $x$, let 
$$\Gamma_x=\Gamma
 {\begin{pmatrix}
   1&x\\0&1
  \end{pmatrix}
}.$$ Then the simplest form of Dani correspondence says that {\it the orbit $\{\Gamma_x a_t\}_{t\geq0}$ is bounded (relatively compact) in $\Gamma\backslash\SL(2,\R)$ if and only if $x$ is a badly approximable number. On the other hand $\{\Gamma_x a_t\}_{t\geq0}$ is divergent if and only if $x$ is rational}. It is well known (see \cite{HW} for instance) that a real number $x$ is badly approximable if and only if the partial quotients in the classical continued fraction of $x$ are bounded. The same is true for $(a,b)$-continued fraction expansion of $x$ with $a=-1$ and $b=1$ because the $(-1,1)$-continued fraction expansion of $x$ is nothing but the classical continued 
fraction expansion of $x$ with alternating signs. In Remark $3.3$ of \cite{CD}, it was shown that a number is badly approximable if and only if the partial quotients in its $(-\frac{1}{2},\frac{1}{2})$-continued fraction expansion are bounded. By the same reasoning (with the help of Proposition \ref{Cusp 2.1}), the same assumption is true for $(a,b)$-continued fraction as well for $(a,b)\in\mathfrak{S}$. So, in the statement of Dani correspondence, the term badly approximable can be replaced by partial quotients in the $(a,b)$-continued fraction expansion of $x$ being bounded, and the rational numbers can be replaced by $x$ having a finite $(a,b)$-continued fraction expansion. Then one 
may think of Theorem \ref{Cusp 1.1} as en extension of Dani correspondence stated above.
\begin{remark}
Let $$E_0=\{x\in\R:\{\Gamma_x a_t\}_{t\geq0}\hspace*{0.2cm}\text
{visits the cusp with frequency 0}\}$$ and 
$$E_\infty=\{x\in\R:\{\Gamma_x a_t\}_{t\geq0}\hspace*{0.2cm}\text{visits the cusp 
with frequency $1$}\}.$$ Then $E_0$ has Hausdorff dimension $1$, since it contains 
the set of badly approximable numbers and it was shown by Jarnik in \cite{J} that 
the set of badly approximable numbers has Hausdorff dimension $1$. On the other hand, it was shown in \cite{FLWW} that if $[a_0,a_1,...]$ (see below for definition) is the classical continued fraction expansion of $x$, then the set of those $x$ for which $\frac{\displaystyle1}
{\displaystyle N}\displaystyle{\sum\limits_{\disp{1}}^{\disp{N}}\log a_j}\rightarrow
\infty$ as $N\rightarrow\infty$, has Hausdorff dimension $\frac{1}{2}$. Then it follows that $E_{\infty}$ has Hausdorff dimension $\frac{1}{2}$. 
\end{remark}

Above remark ensures that $E_\infty$ is a bigger set than the set of rational numbers as the set of rational numbers has Hausdorff dimension $0$. Now we show that $E_\infty$ contains some very well approximable numbers. A real number $x$ is said to be very well approximable if there exists $\varepsilon>0$, such that $|x-\frac{\displaystyle{p}}{\displaystyle{q}}|<
\frac{\displaystyle{1}}{\displaystyle{q^{2+\varepsilon}}}$ holds for infinitely many $q\in\N$ and $p\in\Z$. Recall that each real number $x$ has a classical continued fraction expansion 
\begin{align*}
 x=a_0+\frac{\displaystyle 1}{\displaystyle{a_1+\frac{\displaystyle1}{a_2+
 \frac{\displaystyle1}{\displaystyle{\ddots}}}}},\hspace{0.2cm}
 (a_0\in\Z,a_j\in \N \hspace{0.2cm}\text{for}\hspace{0.2cm} j\geq1),
\end{align*}
written as $x:=[a_0,a_1,a_2,...]$ with $\frac{\displaystyle{p_j}}{\displaystyle{q_j}}=[a_0,a_1,...,a_j]$ denoting the $j$th convergent (see \cite{HW} for more details). Now construct a real number $x=[a_0,a_1,a_2,...]$ using classical continued fraction with the choice of $a_j$'s as follows. Fix some $\varepsilon>0$, choose $a_0\in\Z$ and $a_1\in\N$ arbitrarily, and inductively choose $a_{j+1}=[q_j^\varepsilon]+1$ for $j\geq1$. Then as $\{q_j\}_{j\geq1}$ is an increasing sequence, $\{a_j\}_{j\geq1}$ is also an increasing sequence and consequently
$\frac{\displaystyle1}{\displaystyle N}\displaystyle{\sum\limits_{\disp{1}}^{\disp{N}}}\log a_j
 \rightarrow\infty$ as $N\rightarrow\infty$. Hence, $x\in E_\infty$. On the 
 other hand, it follows from the construction of $x$, that the sequence of convergents 
 $\frac{\displaystyle{p_j}}{\displaystyle{q_j}}$ satisfy the inequality 
 $|x-\frac{\displaystyle{p_j}}{\displaystyle{q_j}}|<\frac{\displaystyle{1}}
 {\displaystyle{q_j^{2+\varepsilon}}}$ for all $j\geq1$, showing that $x$ is a very 
 well approximable number. Note that $E_\infty$ can not contain all very well approximable 
 numbers as the set of very well approximable numbers has Hausdorff 
 dimension $1$ (\cite{J}) and $E_\infty$ has Hausdorff dimension $\frac{\displaystyle1}
 {\displaystyle2}$.
 
\vspace*{0.4cm}

\begin{center} 
{\bf {\Large Summary of revisions over the previous version:}} 
\end{center}
\begin{enumerate}
 \item We have done some rearrangements and modifications in the introduction which leads to a modification in the exposition of the article. Accordingly, we have chosen a better-suited title and modified the abstract.
 \item In the present version of the article, the set $\mathcal{S}$ of parameters $(a,b)$ which is an important object in Theorem \ref{Cusp 1.1}, now has a more restricted description compared to the one in the previous version.  
\end{enumerate}

\section{(a,b)-continued fractions and geodesic flow}
Following S. Katok and I. Ugarcovicci (\cite{KU4}), for $(a,b)\in \mathcal{P}$,
the $(a,b)$-continued fraction expansion of a real number $x$ can be defined using
a generalized integral part function:
\begin{center}
 $[x]_{a,b}:=\left\{
\begin{array}{lr}
[x-a] \hspace*{0.4cm} \text{if} \hspace*{0.2cm} \text{$x<a$}\\
$0$ \hspace*{1.5cm}   \text{if}\hspace*{0.2cm} \text{$a\leq x < b$}\\
\lceil x-b \rceil \hspace*{0.4cm} \text{if} \hspace*{0.2cm} \text{$x\geq b$},
\end{array}\right.$
\end{center}
where  $\lceil x \rceil := [x] + 1$, $[x]$ being the largest integer $\leq x$.  
For $(a,b)\in\mathcal{P}$, every irrational number $x$ can be expressed uniquely as an infinite $(a,b)$-continued fraction of the form (see \cite{KU4} for details)
\begin{align}\label{Expansion1}
 x=a_0-\frac{\displaystyle 1}{\displaystyle{a_1-\frac{\displaystyle1}{a_2-
 \frac{\displaystyle1}{\displaystyle{\ddots}}}}},\hspace{0.2cm}
 (a_j\in\Z,a_j\neq 0 \hspace{0.2cm}\text{for}\hspace{0.2cm} j\geq1),
\end{align}
which we denote by $x:=[a_0,a_1,a_2,...]_{a,b}$. Here $x_0=x$, $a_0=[x_0]_{a,b}$  and $x_j=-\frac{\displaystyle1}{\displaystyle{x_{j-1}-a_{j-1}}}$, $a_j=[x_j]_{a,b}$ for $j\geq1$ and $a_j$ is called the $j$th partial quotient. The rational number $r_j=\frac{\disp{p_j}}{\disp{q_j}}=a_0-\frac{\displaystyle 1}{\displaystyle{a_1-\frac{\displaystyle1}{a_2- \frac{\displaystyle1}{\displaystyle{\ddots-\frac{\displaystyle1}{\displaystyle{a_j}}}}}}}$ is called the $j$th convergent. The sequence $\{|q_j|\}$ is eventually increasing and $r_j$ converges to $x$. As mentioned earlier,
a particular case of $(a,b)$-continued fraction, viz. the $(-1,1)$-continued fraction (also called the alternating continued fraction) is closely related to the classical continued fraction. If $\{a_j\}_{j\geq0}$ is the sequence of partial quotients in the classical continued fraction expansion of a real number $x$, then $\{(-1)^j a_j\}_{j\geq0}$ is the sequence of partial quotients in the $(-1,1)$-continued fraction expansion of $x$. A similar relation holds between the nearest integer continued fraction (also known as Hurwitz's continued fraction introduced by Hurwitz) expansion
and $(-\frac{\displaystyle1}{\displaystyle2},\frac{\displaystyle1}{\displaystyle2})$-continued fraction expansion of any real number. Note that we write the $(a,b)$-continued fraction expansion of any real number as in (\ref{Expansion1}) using minus sign, while in the case of classical or nearest integer continued fraction expansion it is written using plus sign. The use of minus sign while writing the $(a,b)$-continued fraction expansion, presents some advantages which will be clear when we discuss the coding of geodesics using these continued fractions.

Let $\overline{\R}=\R\cup\{\infty\}$ and $f_{a,b}:\overline{\R}\rightarrow\overline{\R}$ be defined by
\begin{center}
 $f_{a,b}(x):=\left\{
 \begin{array}{lr}
   x+1   \hspace*{0.4cm} \text{if} \hspace*{0.2cm} \text{$x<a$}\\
  -\frac{\displaystyle1}{\displaystyle x} \hspace*{0.7cm}   \text{if}\hspace*{0.2cm} \text{$a\leq x < b$}\\
    x-1  \hspace*{0.4cm} \text{if} \hspace*{0.2cm} \text{$x\geq b$}.
 \end{array}\right.$
\end{center}
Note that $f_{a,b}$ is defined using the standard generators $T(x)=x+1$ and 
$S(x)=-\frac{\disp 1}{\disp x}$ of the modular group $\SL(2,\Z)$, and the continued fraction algorithm described above can be obtained using the first
return map of $f_{a,b}$ to the interval $[a,b)$. 
The main object of study in \cite{KU4} is a two dimensional realization of the natural extension
map $F_{a,b}:{\overline{\R}}^2\backslash\Delta\rightarrow{\overline{\R}}^2\backslash
\Delta$, $\Delta=\{(x,y)\in\overline{\R}^2)|x=y\}$ of $f_{a,b}$, which is defined as follows: 
\begin{center}
$F_{a,b}(x,y):=\left\{
\begin{array}{lr}
(x+1,y+1)  \hspace*{0.4cm} \text{if} \hspace*{0.2cm} \text{$x<a$}\\
(-\frac{\displaystyle1}{\displaystyle x},-\frac{\displaystyle1}{\displaystyle y})
\hspace*{1cm}   \text{if}\hspace*{0.2cm} \text{$a\leq x < b$}\\
(x-1,y-1) \hspace*{0.4cm} \text{if} \hspace*{0.2cm} \text{$x\geq b$}.
\end{array}\right.$
\end{center}
The following theorem is just a restatement of the main result of \cite{KU4} for the restricted
set of parameters $\mathcal{P}$.
\begin{theorem}\label{Attractor1}(\cite{KU4}) There exists a one-dimensional Lebesgue measure zero, uncountable set $\mathcal{E}$ contained 
in $\{(a,b)\in\mathcal{P}:b=a+1\}$, such that for all $(a,b)\in\mathcal{P}\backslash\mathcal{E}$,\\\\
$(i)$ the map $F_{a,b}$ has an attractor $D_{a,b}=\bigcap_{n=0}^{\infty}F_{a,b}^{n}(\overline{\R}^2\backslash\Delta)$ on which $F_{a,b}$ is essentially bijective.\\\\
$(ii)$ The set $D_{a,b}$ consists of two (or one in degenerate cases) connected components each having finite rectangular structure, i.e., bounded by non-decreasing step-functions with finitely many steps.\\\\ 
$(iii)$ Every point $(x,y)$ of the plane ($x\neq y$) is mapped to $D_{a,b}$ after finitely many iterations of $F_{a,b}$.
\end{theorem}
In \cite{KU4}, to deduce the above theorem, a crucial role in the arguments used, is played by the orbits of $a$ and $b$ under $f_{a,b}$, viz. to $a$, the upper orbit $\mathcal{O}_u(a)$ (i.e., the orbit of $Sa$) and the lower orbit $\mathcal{O}_l(a)$ (i.e., the orbit of $Ta$), and to $b$, the upper orbit $\mathcal{O}_u(b)$ (i.e., the orbit of $T^{-1}b$) and the lower orbit $\mathcal{O}_l(b)$ (i.e., the orbit of $Sb$). Let us denote the set $\mathcal{P}\backslash\mathcal{E}$ by the symbol $\mathcal{S}$. It was proved in \cite{KU4} that if $(a,b)\in\mathcal{S}$, then $f_{a,b}$ satisfies the finiteness condition. This means that for both $a$ and $b$, their upper and lower orbits are either eventually periodic, or they satisfy the cycle property, i.e., they meet forming a cycle, in other words there exist integers $k_1,m_1,k_2,m_2\geq0$ such that 
$$f_{a,b}^{m_1}(Sa)=f_{a,b}^{k_1}(Ta)=c_a \hspace*{0.1cm} (\text{respectively}\hspace{0.2cm}f_{a,b}^{m_2}(T^{-1}b)=f_{a,b}^{k_2}(Sb)=c_b),$$ where $c_a$ and $c_b$ are the ends of the cycles. If the products of transformations over the upper and lower sides of the cycle of $a$ (respectively $b$) are equal, $a$ (respectively $b$) is said to have strong cycle property, otherwise it has weak cycle property. Let 
\begin{center}
$\mathcal{L}_a=\left\{
\begin{array}{lr}
\mathcal{O}_l(a) \hspace*{4.5cm} \text{if $a$ has no cycle property}\\
\text{lower part of $a$-cycle}\hspace*{1.8cm} \text{if $a$ has strong cycle property}\\
\text{lower part of $a$-cycle}\cup\{0\}\hspace*{0.7cm} \text{if $a$ has weak 
cycle property},
\end{array}\right.$
\end{center}

\begin{center}
$\mathcal{U}_a=\left\{
\begin{array}{lr}
\mathcal{O}_u(a) \hspace*{4.5cm} \text{if $a$ has no cycle property}\\
\text{upper part of $a$-cycle}\hspace*{1.8cm} \text{if $a$ has strong cycle property}\\
\text{upper part of $a$-cycle}\cup\{0\}\hspace*{0.7cm} \text{if $a$ has weak 
cycle property}
\end{array}\right.$
\end{center}
and $\mathcal{L}_b$, $\mathcal{U}_b$ be defined similarly. Also let $\mathcal{L}_{a,b}=
\mathcal{L}_a\cup\mathcal{L}_b$ and $\mathcal{U}_{a,b}=\mathcal{U}_a\cup\mathcal{U}_b$. So, $f_{a,b}$ satisfies the finiteness condition means that both the sets $\mathcal{L}_{a,b}$ and $\mathcal{U}_{a,b}$ are finite, which is true when  $(a,b)\in\mathcal{S}$.
In \cite{KU4}, first a set $A_{a,b}$, having finite rectangular structure, was constructed (see Theorem $5.5$ in \cite{KU4}) using the values in the sets $\mathcal{U}_{a,b}$ and $\mathcal{L}_{a,b}$, and then it was shown (Theorem $6.4$ in \cite{KU4}) that $A_{a,b}$ actually coincides with the attractor $D_{a,b}$. The upper component of $D_{a,b}$ is bounded by non-decreasing step functions with values in the set 
$\mathcal{U}_{a,b}$ and the lower component of $D_{a,b}$ is bounded by non-decreasing step functions with values in the set $\mathcal{L}_{a,b}$. 

Making use of the properties of the map $F_{a,b}$ and the attractor $D_{a,b}$, in a subsequent paper (\cite{KU5}), S. Katok and I. Ugarcovicci developed a general 
method of coding geodesics on the modular surface and gave a symbolic description of the geodesic flow associated with the modular surface. We first recall from 
\cite{KU5}, the notion of $(a,b)$-reduced geodesics, which plays a crucial role in determining the cross-section for the geodesic flow needed for coding purposes.
\begin{definition}
 A geodesic in $\H$ with real endpoints $u$ and $w$, $w$ being the attracting and $u$
 being the repelling endpoints, is called $(a,b)$-reduced if $(u,w)\in\Lambda_{a,b}$,
 where 
 \begin{align*}
\Lambda_{a,b}:=F_{a,b}(D_{a,b}\cap\{a\leq w < b\})=S(D_{a,b}\cap\{a\leq w < b\}).
\end{align*}
\end{definition}
Given any geodesic $\gamma'$ in $\H$, one can obtain an $(a,b)$-reduced geodesic $\PSL(2,\Z)$-equivalent to $\gamma'$ by using the reduction property ($3$rd assertion in Theorem \ref{Attractor1}) of the map $F_{a,b}$. More precisely, if $\gamma'$ is a geodesic which is not $(a,b)$-reduced and if $w'=[a'_0,a'_1,a'_2,...]_{a,b}$ is the attracting end point of $\gamma'$, then there exists some positive integer $n$ such that $ST^{-a'_n}...ST^{-a'_1}ST^{-a'_0}(\gamma')$ is an $(a,b)$-reduced geodesic (see \cite{KU5} for details). Now let $\gamma$ be an $(a,b)$-reduced geodesic with attracting and repelling endpoint $w$ and $u$ respectively, and $[a_0,a_1,a_2,...]_{a,b}$ be the $(a,b)$-continued fraction expansion of $w$. Using the essential bijectivity of the map $F_{a,b}$, one can extend the sequence $(a_0,a_1,a_2,...)$ in the past as well to get a bi-infinite sequence $(...,a_{-2},a_{-1},a_0,a_1,a_2,...)$, called the coding sequence of $\gamma$ and written as 
\begin{center}
 $[\gamma]_{a,b}=(...,a_{-2},a_{-1},a_0,a_1,a_2,...)$,
\end{center}
where $a_{-1}-\frac{\displaystyle 1}{\displaystyle{a_{-2}-
\frac{\displaystyle1}{\displaystyle{\ddots}}}}=\frac{\displaystyle1}{\displaystyle u}$ 
(see Section $3$ of \cite{KU5} for details).

Now we recall from \cite{KU5}, the description of the cross-section. Let
\begin{center}
 $C=\{z\in\H|\hspace*{0.2cm} |z|=1,\text{Im}z\geq0\}$
\end{center}
be the upper half of the unit circle and $\mathcal{F}$ denote the standard fundamental domain for the action of $\SL(2,\Z)$ on $\H$, given by 
\begin{center}
$\mathcal{F}:=\{z=x+iy\in\H|\hspace*{0.2cm}|z|\geq1, |x|\leq \frac{\disp1}{\disp2}\}$.
\end{center}
\begin{figure}\label{pic4}
\centering
\thispagestyle{empty}
\includegraphics[scale=0.7]{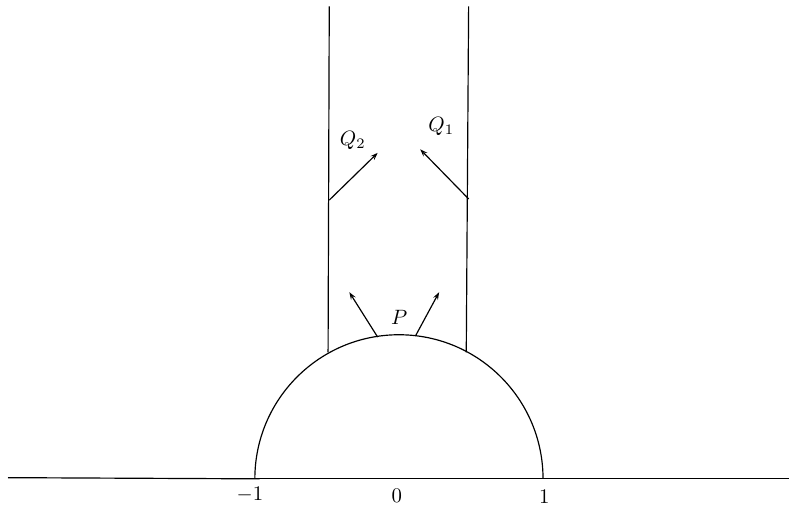}
\caption{Cross-section $C_{a,b}$}
\end{figure}
Using the definition of $(a,b)$-reduced geodesic it is easy see the following fact.
\begin{proposition}\label{Coding1}(\cite{KU5})  For $(a,b)\in\mathcal{S}$, every $(a,b)$-reduced geodesic intersects $C$.
\end{proposition}
Given an $(a,b)$-reduced geodesic $\gamma$ with attracting and repelling endpoints $w$ and $u$ respectively, the {\it cross-section point} on $\gamma$ is the 
intersection point of $\gamma$ with $C$. Let $\phi :\Lambda_{a,b}\rightarrow T^1\H$ be defined by 
\begin{center}
 $\phi(u,w):=(z,\zeta)$,
\end{center}
where $z\in\H$ is the cross-section point on the geodesic $\gamma$ joining $u$ and $w$, and $\zeta$ is the unit vector tangent to $\gamma$ at $z$. The map $\phi$ is clearly injective and after composing with the Canonical projection $\pi$ we obtain a map
\begin{center}
 $\pi\circ\phi:\Lambda_{a,b}\rightarrow T^1M$.
\end{center}
Let $C_{a,b}:=\pi\circ\phi(\Lambda_{a,b})\subset T^1M$. Then $C_{a,b}$ is a cross-section for the geodesic flow associated with the modular surface (see \cite{KU5} for details). The lift of $C_{a,b}$ to $T^1\H$ restricted to the unit tangent vectors having base points on the fundamental domain $\mathcal{F}$, can be described as follows:
\begin{center}
$\pi^{-1}(C_{a,b})\cap\left(\bigcup\limits_{z\in\mathcal{F}}T_z^1\H\right)=P\cup Q_1 \cup Q_2$ (see Figure $1$),
\end{center}
where $P$ consists of unit tangent vectors on the circular boundary of the fundamental region $\mathcal{F}$ and pointing inward such that the corresponding geodesic $\gamma$
on $\H$ is $(a,b)$-reduced; $Q_1$ consists of unit tangent vectors with base points on the right vertical boundary of $\mathcal{F}$ and pointing inward such that if $\gamma$ is the geodesic corresponding to one such unit vector, then $TS\gamma$ is $(a,b)$-reduced; $Q_2$ consists of unit tangent vectors with base points on the left vertical boundary of $\mathcal{F}$ and pointing inward such that if 
$\gamma$ is the geodesic corresponding to one such unit vector, then $T^{-1}S\gamma$ is $(a,b)$-reduced. 

Now let $v\in T^1M$ and $\gamma_v$ be the corresponding geodesic in $M$ and 
$\tilde{\gamma}_v$ be an $(a,b)$-reduced lift of it inside $\H$. Also let $\eta:T^1M\rightarrow M$ be the Canonical projection of $T^1M$ onto $M$. The following theorem from \cite{KU5} provides the base for coding geodesics on the modular surface using $(a,b)$-continued fractions. 
\begin{theorem}\label{Coding2}(\cite{KU5})
Let $\gamma_v$ and $\tilde{\gamma}_v$ be as above. Then each geodesic segment of $\gamma_v$ between successive returns to 
$\eta(C_{a,b})$, while extended to a geodesic, produces an $(a,b)$-reduced geodesic on $\H$, and each $(a,b)$-reduced geodesic $\PSL(2,\Z)$-equivalent to $\tilde{\gamma}_v$ is obtained 
in this way.The first return of $\gamma_v$ to $\eta(C_{a,b})$ corresponds to a left shift of the coding sequence of $\tilde{\gamma}_v$.
 \end{theorem}
Let $\{g_tv\}$ be the orbit of the geodesic flow on $T^1M$ corresponding to the geodesic $\gamma_v$, i.e., $\gamma_v(t)=\eta(g_tv)$ and let $\gamma_v^j$ be the segment of the geodesic $\gamma_v$ corresponding to the portion of the orbit $\{g_tv\}_{t\geq0}$ between $(j-1)$th and $j$th returns to the cross-section $C_{a,b}$. We call the segment $\gamma_v^j$ the $j$th {\it excursion} of the geodesic $\gamma_v$ into the cusp. Let $w=[a_0,a_1,a_2,...]_{a,b}$ be the attracting end point of $\tilde{\gamma}_v$ and $\tilde{\gamma}_{vj} := ST^{-a_{j-1}}...ST^{-a_1}ST^{-a_0}(\tilde{\gamma}_v)$. Then the segment of $\tilde{\gamma}_{vj}$ between $C$ and $a_j+C$, denoted by $\tilde{\gamma}_v^j$ is a lift of $\gamma_v^j$ to the hyperbolic plane. Assuming the geodesics to be parameterized by arc length, the time between the $(j-1)$th and the $j$th return of $\{g_tv\}$ to the cross-section, called the $j$th {\it return time}, is given by $$t_j:=h(\gamma_v^j)=h(\tilde{\gamma}_v^j),$$ where $h$ stands for the hyperbolic length of the geodesic segment. Also let $\mathcal{S}':=\mathcal{S}\backslash (-1,1)$.
\begin{proposition}\label{Cusp 2.1}
 If $(a,b)\in \mathcal{S}'$, then $C_{a,b}$ is contained inside a compact subset of $T^1M$.
\end{proposition}
\begin{proof}
 The structure of $D_{a,b}$ is discussed in detail in Theorem $5.5$ of \cite{KU4}. $D_{a,b}$ has two connected components, the lower one we denote by $D_{a,b}^l$ and the upper one we
 denote by $D_{a,b}^u$. Both the sets $D_{a,b}^l$ and $D_{a,b}^u$ have finite rectangular structure i.e., bounded by non-decreasing step functions with finite number of steps. For $D_{a,b}^l$ the values of the step function are given by the set $\mathcal{L}_{a,b}$, and for $D_{a,b}^u$ the values of the step function are given by the set $\mathcal{U}_{a,b}$. The structure of the boundary (see Figure $2$ for a typical picture of $D_{a,b}$) of $D_{a,b}$ consists
\begin{figure}\label{pcusp1}
\centering
\thispagestyle{empty}
\includegraphics[scale=0.7]{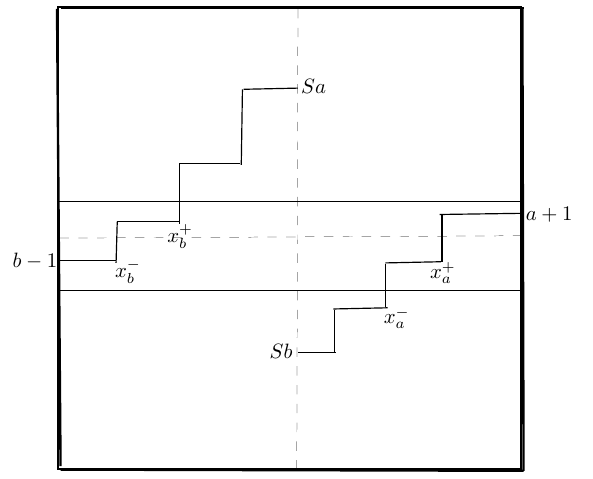}
\caption{Structure of $D_{a,b}$}
\end{figure}
 of finite number of horizontal segments at different points of the set $\mathcal{L} _{a,b}$, called the different levels of $D_{a,b}^l$ and consecutive levels are joined by vertical segments, where the highest level is $y=a+1$. $D_{a,b}^u$ has a similar description with the lowest level being $y=b-1$. Let $x_a^-$ be the $x$-coordinate of the vertical segment joining two consecutive levels $y_a^-$ and $y_a^+$ of $D_{a,b}^l$ with $y_a^-\leq a < y_a^+$, and $x_a^+$ be the $x$-coordinate of the vertical segment joining two consecutive levels $y_-$ and $y_+$ with $y_-\leq 0 < y_+$. Similarly, let $x_b^-$ be the $x$-coordinate of the vertical segment joining two consecutive levels $y_-'$ and $y_+'$ of $D_{a,b}^u$  with $y_-'<0\leq y_+'$, and $x_b^+$ be the $x$-coordinate of the vertical segment joining two consecutive levels $y_b^-$ and $y_b^+$ with $y_b^-<b\leq y_b^+$. Also let $y_l$ be the level above $Sb$ and next to $Sb$; $y_u$ be the level below $Sa$ and next to $Sa$.
 
 It follows from these assertions and the definition of $\Lambda_{a,b}$, that a geodesic $\tilde{\gamma}_v$ with attracting and repelling endpoints $w$ and $u$ respectively with $w>0$, is $(a,b)$-reduced if and only if 
 $$(u,w)\in\left[-\frac{\displaystyle1}{\displaystyle{x_a^-}},0\right)\times 
  \left[-\frac{\displaystyle1}{\displaystyle a},\infty\right)\bigcup
  \left(0,-\frac{\displaystyle1}{\displaystyle{x_b^-}}\right]\times 
  \left[-\frac{\displaystyle1}{\displaystyle{b-1}},\infty\right).$$
 On the other hand if $w<0$, then $\tilde{\gamma}_v$ is $(a,b)$-reduced if and only if
 $$(u,w)\in\left(0,-\frac{\displaystyle1}{\displaystyle{x_b^+}}\right]\times 
\left(-\infty,-\frac{\displaystyle1}{\displaystyle b}\right]\bigcup
\left[-\frac{\displaystyle1}{\displaystyle{x_a^+}},0\right)\times 
\left(-\infty,-\frac{\displaystyle1}{\displaystyle{a+1}}\right].$$
We show that $x_a^-$, $x_a^+>1$ and $x_b^-$, $x_b^+<-1$.
  
For $(a,b)\in \mathcal{S}'$, let $m_a$ and $m_b$ be positive
integers such that $a\leq T^{m_a}STa<a+1$ and $a\leq T^{m_b}Sb<a+1$. Let $m_a$, $m_b\geq3$,
then the the proof of Lemma $5.6$ of \cite{KU4} shows that the vertical segment joining $Sb$ and $y_l$ has $x$-coordinate greater than $1$, and the vertical segment
joining $y_u$ and $Sa$ has $x$-coordinate less than $-1$. Therefore, in these cases we have
$x_a^-$, $x_a^+>1$ and $x_b^-$, $x_b^+<-1$. Now we consider the situation when $m_a$, $m_b\leq2$. Note that $m_a$ can never be $1$, for if $m_a=1$, then $a=0$ since $a>-1$, but we have assumed that $a<0$. So, $m_a\geq2$. Now if either $m_a$ or $m_b$ is $2$, then from the explicit cycle description of $a$ and $b$ discussed in \cite{KU4}, we see that there is always one level between $y_l$ and $a$; similarly there is always one level between $b$ and $y_u$. As the statement of Lemma $5.6$ of \cite{KU4}
guarantees that the vertical segment joining $Sb$ and $y_l$ has $x$-coordinate greater than or 
equal to $1$ and the vertical segment joining $y_u$ and $Sa$ has $x$-coordinate less than or 
equal to $-1$, it follows that $x_a^-$, $x_a^+>1$ and $x_b^-$, $x_b^+<-1$ in these cases as 
well.  
 \begin{figure}\label{hd2}
\centering
\thispagestyle{empty}
\includegraphics[scale=1]{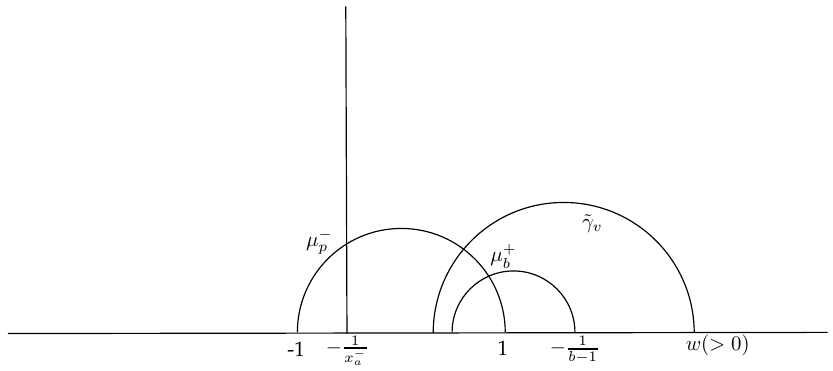}
 \caption{Cross-section point of an $(a,b)$-reduced geodesic}
\end{figure} 

From the discussion above we have, $-\frac{\displaystyle{1}}{\displaystyle{x_b^-}}<1$ and 
$-\frac{\displaystyle{1}}{\displaystyle{x_a^-}}>-1$.
Now let $\mu_a^+$ be the intersection point of the geodesic joining $0$ and $-\frac{\displaystyle1}{\displaystyle{a}}$, and $C$; $\mu_b^+$ be the intersection point of the geodesic joining $-\frac{\displaystyle{1}}{\displaystyle{x_b^-}}$ and 
$-\frac{\displaystyle{1}}{\displaystyle{b-1}}$, and $C$.
We chose one of $\mu_a^+$ and $\mu_b^+$, which has $y$-coordinate less than or equal to the other and denote it by $\mu_p^+$. Also let $\mu_p^-$ be the intersection point of $C$ and the vertical geodesic based at the point $-\frac{\displaystyle{1}}
{\displaystyle{x_a^-}}$. Then any $(a,b)$-reduced geodesic $\tilde{\gamma}_v$ having attracting endpoint $w>0$, intersects the segment joining $\mu_p^-$ and $\mu_p^+$ of $C$ (see Figure $3$). Consequently the cross-section point for any $(a,b)$-reduced geodesic having positive attracting endpoint, 
has $y$-coordinate uniformly bounded away from $0$.  The same is true for any $(a,b)$-reduced geodesic with negative attracting endpoint as well, which can be shown similarly by using the fact that 
$-\frac{\displaystyle{1}}{\displaystyle{x_b^+}}<1$ and 
$-\frac{\displaystyle{1}}{\displaystyle{x_a^+}}>-1$. This completes the proof of the proposition.
\end{proof}

\section{Cusp excursions with extreme frequencies}
In this section, we prove the main results of this article. The results are about classifying two kinds of forward orbits of geodesic flow apart from the generic ones. This is done by relating the time spent by the orbits in cusp neighbourhoods compared to the total time parameter, and the average growth rate of the partial quotients of the continued fraction expansion of the attracting end points of the corresponding geodesics. It is worth mentioning that there are many interesting results relating cusp excursions of geodesics on hyperbolic 2-orbifolds and Diophantine approximation. For example, see \cite{H3}, \cite{H4} and the references given there (As there is a large body of literature around this phenomena, the reference list given here is not complete by any means). In \cite{H3} and \cite{H4}, various aspects of cusp excursions of a generic set of geodesics have been studied and analogue of various results from classical Diophantine approximation in the context of Fuchsian groups have been obtained, while restricting to the case of the modular surface these produce new proofs of classical results (see \cite{H3}, \cite{H4} for details). For example, it was shown in \cite{H4} that for $d<1$ and for almost all $v\in T^1M$, if $\{\gamma_v^{j_k}\}$ is the subsequence of $\{\gamma_v^j\}$ which intersect $M_d$, then $\lim_{n\rightarrow\infty}\frac{1}{n}\sum\limits_{k=1}^n h(M_d\cap \gamma_v^{j_k})=\pi$.
 
In this article, we consider a certain class of geodesics apart from the generic ones and look at their behaviour
in terms of spending time inside cusp neighbourhoods compared to their length parameter. Let $v$, $\gamma_v$, $\tilde{\gamma}_v$,
$w=[a_0,a_1,a_2,...]_{a,b}$ be as in the previous section. The partial quotients of the continued fraction expansion of $w$ determine how much further the orbit $\{g_tv\}$ of the geodesic flow goes into a typical neighbourhood of the cusp before returning to the cross-section. This particular fact is easier to see
when the cross-section  $C_{a,b}$ is contained inside a compact set which is the case when $(a,b)\in\mathcal{S}'$.
Whereas for $(a,b)=(-1,1)$, the cross-section $C_{-1,1}$ is not contained inside a compact set. In this case 
we use the formula for return times given by S. Katok and I. Ugarcovicci and some other facts which are particular 
to the $(-1,1)$-continued fraction.

\subsection{$(a,b)\in\mathcal{S}'$}
\begin{lemma}\label{Cusp 3.1}
Let $j$ be a positive integer.\\\\ $(i)$ Assume that $a_j>0$. Then $\tilde{\gamma}_{vj}$ intersects or does not intersect $\H_d$ accordingly as $a_j>2d-a-\frac{\displaystyle{1}} {\displaystyle{x_b^-}}$ or  
$a_j<2d-b-\frac{\displaystyle{1}}
{\displaystyle{x_a^-}}$.\\\\
$(ii)$ Assume that $a_j<0$. Then $\tilde{\gamma}_{vj}$ intersects or does not intersect $\H_d$ accordingly as $|a_j|>2d+b+\frac{\displaystyle{1}}{\displaystyle{x_a^+}}$ or $|a_j|<2d+a+\frac{\displaystyle{1}}{\displaystyle{x_b^+}}$.  
\end{lemma}
\begin{proof}
If $a_j>0$, then the attracting endpoint $w_j$ of $\tilde{\gamma}_{vj}$ lies in the interval\\ $[a_j+a,a_j+b)$ and the repelling endpoint $u_j$ is contained in the interval $\left(-\frac{\displaystyle{1}}{\displaystyle{x_a^-}},-\frac{\displaystyle{1}}
{\displaystyle{x_b^-}}\right)$. So, in this case, $\tilde{\gamma}_{vj}$ lies above the geodesic  $\tilde{\gamma}_{a_j}^-$, where $\tilde{\gamma}_{a_j}^-$ is the geodesic joining $-\frac{\displaystyle{1}}{\displaystyle{x_b^-}}$ and $a_j+a$. Also, $\tilde{\gamma}_{vj}$ lies  below the geodesic $\tilde{\gamma}_{a_j}^+$, where $\tilde{\gamma}_{a_j}^+$ is the geodesic joining 
$-\frac{\displaystyle{1}}{\displaystyle{x_a^-}}$ and $a_j+b$. So, if the radius of $\tilde{\gamma}_{a_j}^+$ is less than $d$, then $\tilde{\gamma}_{vj}$ does not intersect $\H_d$; on the other hand if the radius of $\tilde{\gamma}_{a_j}^-$ is greater than $d$, then $\tilde{\gamma}_{vj}$ does intersect $\H_d$. Now a simple calculation gives the assertion of the lemma in the case $a_j>0$. If $a_j<0$, then $\tilde{\gamma}_{vj}$ lies above the geodesic joining $a_j+b$ and $-\frac{\displaystyle{1}}
{\displaystyle{x_a^+}}$, and lies below the geodesic joining $a_j+a$ and $-\frac{\displaystyle{1}}{\displaystyle{x_b^+}}$. Again a simple calculation gives the assertion of the lemma in the case $a_j<0$.
\end{proof}
The following two lemmas which are crucial to the arguments to follow, can be proved easily using the fact that the cross-section $C_{a,b}$ is contained inside a compact set in $T^1M$. The proof of similar statements for the particular case $(a,b)=\left(-\frac{\displaystyle1}{\displaystyle2},
\frac{\displaystyle1}{\displaystyle2}\right)$ is contained in \cite{CD} (Proposition $3.4$ and Proposition $3.5$ respectively) and the same proofs work for any $(a,b)\in\mathcal{S}'$ as well.  
\begin{lemma}\label{Cusp 3.2}
Let $d>1$ be such that $\overline{M}_d\cap C_{a,b}=\emptyset$, then if 
$\gamma_v^j\cap M_d$ is nonempty, $\tilde{\gamma}_v^j\cap\H_d$ is the only 
connected component of $\pi^{-1}(\gamma_v^j\cap M_d)$.
\end{lemma}

\begin{lemma}\label{Cusp 3.3}
 Let $v\in T^1M$, $\gamma_v$ be the corresponding geodesic in $M$, and $\tilde{\gamma}_v$ be an $(a,b)$-reduced lift of $\gamma_v$ inside $\H$. Let $w=[a_0,a_1,a_2,...]_{a,b}$ be the attracting end point of $\tilde{\gamma}_v$ and $t_j$ be the $j$th return time for the corresponding orbit $\{g_tv\}$ of the geodesic flow. Then there exist a constant $\kappa>0$ such that 
 \begin{align}\label{equation2}
|t_j-2\log |a_j||\leq \kappa, \hspace*{0.1cm}\forall \hspace*{0.1cm} j\geq 0.
\end{align}
\end{lemma}

\begin{remark}
The asymptotic estimates for values of binary quadratic forms at integer points were obtained in \cite{CD} in terms of $(-\frac{\displaystyle1}{\displaystyle2},
\frac{\displaystyle1}{\displaystyle2})$-continued fraction expansion of the coefficients of the quadratic forms, and the $(-\frac{\displaystyle1}{\displaystyle2}, \frac{\displaystyle1}{\displaystyle2})$-continued fraction coding of geodesics on the modular surface was used to obtain the estimates. The facts that the cross-section for geodesic flow corresponding to the $(-\frac{\displaystyle1}{\displaystyle2},
\frac{\displaystyle1}{\displaystyle2})$-continued fraction coding, is contained inside a compact subset of $T^1M$ and the return times can be bounded uniformly by the partial quotients as in (\ref{equation2}), were used crucially to obtain those estimates. Since the above two properties hold for $(a,b)$-continued fraction coding as well for $(a,b)\in\mathcal{S}'$, one can obtain similar estimates as in \cite{CD} for values of binary quadratic forms at integer points in terms of the $(a,b)$-continued fraction expansions of its coefficients as well.
\end{remark}
Given $d>1$, let $\underline{d}^+=2d-b-\frac{\displaystyle{1}}{\displaystyle{x_a^-}}$,
 $\underline{d}^-=2d+a+\frac{\displaystyle{1}}{\displaystyle{x_b^+}}$, and 
\begin{center}
$\mathfrak{j}_{\underline{d}}^N=\#\big(0\leq j<N:\text{either}\hspace*{0.1cm}a_j>
\underline{d}^+ \hspace*{0.2cm}\text{if}\hspace*{0.2cm}a_j>0\hspace*{0.2cm}
\text{or}\hspace*{0.2cm}a_j<-\underline{d}^-
\hspace*{0.2cm}\text{if}\hspace*{0.2cm}a_j<0\big)$.
\end{center}
Let $\bar{d}_+=2d-a-\frac{\displaystyle{1}}{\displaystyle{x_b^-}}$, 
$\bar{d}_-=2d+b+\frac{\displaystyle{1}}{\displaystyle{x_a^+}}$, and
\begin{center}
$\mathfrak{j}_{\bar{d}}^N=\#\big(0\leq j<N:\text{either}\hspace*{0.1cm}a_j>
\bar{d}_+ 
\hspace*{0.2cm}\text{if}\hspace*{0.2cm}a_j>0\hspace*{0.2cm}\text{or}\hspace*{0.2cm}a_j
<-\bar{d}_-\hspace*{0.2cm}\text{if}\hspace*{0.2cm}a_j<0\big)$.
\end{center}
Let $S_N=t_1+t_2+...+t_N$. Also let $$I_N^d:=\frac{\disp{1}}{\disp{S_N}}\disp{\int_{\disp{0}}^{\disp{S_N}}\chi_d(g_tv)dt}, \hspace*{0.3cm}I_T^d:=\frac{\disp{1}}{\disp{T}}\disp{\int_{\disp{0}}^{\disp{T}}\chi_d(g_tv)dt},$$ where $\chi_d$ denotes the characteristic function of the neighbourhood $\overline{M}_d$ of the cusp and $v\in T^1M$.

It is evident from Lemma \ref{Cusp 3.1}, Lemma \ref{Cusp 3.2} and Lemma \ref{Cusp 3.3} that the $j$th excursion of the geodesic goes more and more into the cusp as the value of $|a_j|$ gets bigger and bigger and vice versa. The following proposition uses this fact to characterize those orbits of geodesic flow which visit the cusp with full frequency. It is easy to see that to conclude about the extreme behaviour of $I_T^d$, it is enough to consider $I_N^d$.
\begin{proposition}\label{Cusp 3.5}
Let $v\in T^1M$, $\gamma_v$ be the corresponding geodesic on $M$ and $\tilde{\gamma}_v$ be an $(a,b)$-reduced lift of $\gamma_v$ in $\H$. Let $w=[a_0,a_1,a_2,...]_{a,b}$ be the attracting endpoint of $\tilde{\gamma}_v$. Then
$I_N^d=\frac{\disp{1}}{\disp{S_N}}\disp{\int_{\disp{0}}^{\disp{S_N}}\chi_d(g_tv)dt}\rightarrow 
1$ as $N\rightarrow\infty$ for all $d>1$, if and only if,\\
$\frac{\displaystyle 1}{\displaystyle N}(\log|a_0|+\log|a_2|.......+\log|a_{N-1}|) \rightarrow\infty$ as $N\rightarrow\infty$.
\end{proposition}
\begin{proof}
We enumerate those $j$ for which either $a_j>\bar{d}_+$ for $a_j>0$, or $a_j<-\bar{d}_-$ for \\ $a_j<0$, by the subsequence $\{\disp{j_k}\}$, and by
$\sum\limits_{\disp{k=1}}^{\disp{\mathfrak{j}_{\bar{d}}^N}}
\log|a_{\disp{\disp{j_k}}}|$ we mean the sum
\begin{center}
$\disp{\sum\limits_{\disp{a_j>\bar{d}_+
\hspace*{0.1cm}\text{or}\hspace*{0.1cm}a_j<-\bar{d}_-,0\leq j\leq N-1}}\log|a_j|}$.
\end{center}
On the other hand, we enumerate those $j$ for which $a_j\leq\bar{d}_+$ if $a_j>0$, or $a_j\geq-\bar{d}_- $ if $a_j<0$, by the subsequence $\{j_l\}$, and by
$\sum\limits_{\disp{l=1}}^{\disp{N-\mathfrak{j}_{\bar{d}}^N}}\log|a_{\disp{j_l}}|$ we mean the sum 
\begin{center}
$\disp{\sum\limits_{\disp{0<a_j\leq\bar{d}_+
\hspace*{0.1cm}\text{or}\hspace*{0.1cm}-\bar{d}_-\leq a_j<0,0\leq j\leq N-1}}\log|a_j|}$.
\end{center}
Now suppose
$\frac{\disp 1}{\disp N}\disp{\sum\limits_{\disp{{j=0}}}^{\disp{N-1}}\log{|a_j|}}
\rightarrow\infty$ as $N\rightarrow\infty$ which implies by Lemma \ref{Cusp 3.3}, that
$\frac{\disp1}{\disp N}\disp{\sum\limits_{\disp{j=1}}^{\disp N} t_j} \rightarrow\infty$ as $N\rightarrow\infty$.
 
Let $d_{a,b}>0$ be such that $\overline{M}_{\disp{d_{a,b}}}\cap
C_{a,b}=\emptyset$. Now for any $d>d_{a,b}$, let 
$$c_{\disp{\disp{j_k}}}:=h(M_d\backslash\gamma_v^{j_k}).$$ Then
\begin{align}
I_N^d=\frac{\disp1}{\disp{S_N}}\int_{\disp{0}}^{\disp{S_N}}\chi_d(g_tv)dt
&\geq \frac{\disp{\sum\limits_{\disp{k=1}}^{\disp{\mathfrak{j}_{\bar{d}}^N}}}
\disp{(t_{\disp{\disp{j_k}}}-c_{\disp{\disp{j_k}}})}}{\disp{\sum\limits_{\disp{j=1}}^{\disp{N}} t_j}}\\
&=1-\frac{\disp{\frac{1}{N}\sum\limits_{\disp{l=1}}^{\disp{N-\mathfrak{j}_{\bar{d}}^N}}}
\disp{t_{\disp{\disp{j_l}}}}}{\disp{\frac{1}{N}\sum\limits_{\disp{j=1}}^{\disp{N}} t_j}}-
\frac{\disp{\frac{1}{N}\sum\limits_{\disp{k=1}}^{\disp{\mathfrak{j}_{\bar{d}}^N}}}
\disp{c_{\disp{\disp{j_k}}}}}{\disp{\frac{1}{N}\sum\limits_{\disp{j=1}}^{\disp{N}} t_j}},
\end{align}
As both the quantities
$\frac{\disp1}{\disp N}\disp{\sum\limits_{\disp{l=1}}^{\disp{N-\mathfrak{j}_{\bar{d}}^N}}t_{\disp{\disp{j_l}}}}$ and
$\frac{\disp1}{\disp N}\disp{\sum\limits_{\disp{k=1}}^{\disp{\mathfrak{j}_{\bar{d}}^N}}}
\disp{c_{\disp{\disp{j_k}}}}$ are bounded, and\\ $\frac{\disp1}{\disp N}\disp{\sum\limits_{\disp{j=1}}^{\disp{N}} t_j}\longrightarrow\infty$ as $N\rightarrow\infty$, it follows that $I_N^d\rightarrow 1$ as $N\rightarrow\infty$.

To prove the converse statement, we show that if
$\frac{\disp1}{\disp N}\disp{\sum\limits_{\disp{{j=0}}}^{\disp{N-1}}\log{|a_{\disp{j}}|}}\nrightarrow
\infty$ as $N\rightarrow\infty$, then there is some $\d>1$ such that $I_N^d$ can not go to $1$ as $N\rightarrow\infty$. Now  $\frac{\disp1}{\disp N}\disp{\sum\limits_{\disp{{j=1}}}^
{\disp{N}}\log{|a_j|}}\nrightarrow\infty$ as $N\rightarrow\infty$ means that there is a subsequence $\{\disp{N_s}\}$ and $\mathfrak{m}>0$ such that
$\frac{\disp1}{\disp{\disp{N_s}}}\disp{\sum\limits_{\disp{j=1}}^{\disp{\disp{N_s-1}}}\log|a_j|}
<\mathfrak{m}$ for all $s\in\mathbb{N}$, which again means, by Lemma \ref{Cusp 3.3}, that $\frac{\disp1}{\disp{\disp{N_s}}}\disp{\sum\limits_{\disp{j=1}}^{\disp{\disp{N_s}}}
t_j}<\tilde{\mathfrak{m}}$ for some $\tilde{\mathfrak{m}}>0$ and for all $s\in\mathbb{N}$. Since 
$\frac{\disp1}{\disp{N}}\mathfrak{j}_{\disp{\bar{d}}}^{\disp{N}}\rightarrow 1$ as $N\rightarrow\infty$ for all $d>1$ implies $\frac{\disp1}{\disp N}\disp{\sum\limits_{\disp{{j=1}}}
^{\disp{N}} t_j}\rightarrow\infty$ as $N\rightarrow\infty$, which again by Lemma \ref{Cusp 3.3} implies 
$\frac{\disp1}{\disp N}\disp{\sum\limits_{\disp{{j=0}}}^{\disp{N-1}}
\log{|a_j|}}\rightarrow\infty$ as $N\rightarrow\infty$,
we may assume that there exists some $r>0$ and $d>1$ such that (if needed by considering a subsequence of $\{N_s\}$ and denoting it again by 
$\{N_s\}$) $\frac{\disp1}{\disp{\disp{N_s}}}\mathfrak{j}_{\disp{\bar{d}}}^{\disp{\disp{N_s}}}
 <1-r$ for all $s\in\mathbb{N}$.
Now
\begin{align}
I_{\disp{\disp{N_s}}}^d\leq 1-\frac{\disp{\frac{1}{\disp{N_s}}\sum\limits_{\disp{l=1}}^
{\disp{\disp{N_s}-\mathfrak{j}_{\bar{d}}^{\disp{\disp{N_s}}}}}}\disp{t_{\disp{\disp{j_l}}}}}{\disp{\frac{1}{\disp{N_s}}
\sum\limits_{\disp{j=1}}^{\disp{\disp{N_s}}} t_j}}.
\end{align}
Since the cross-section point for any $(a,b)$-reduced geodesic, is uniformly bounded away from the real line, it follows that $t_j$ has a uniform lower bound, i.e., $t_j>\mathfrak{t}$ for some $\mathfrak{t}>0$ and all $j\geq0$. Since $\frac{\disp1}{\disp{\disp{N_s}}}(\disp{N_s}-\mathfrak{j}_{\bar{d}}^
{\disp{\disp{N_s}}})>r$ and $\frac{\disp1}{\disp{N_s}}\disp{\sum\limits_{\disp{j=1}}
^{\disp{N_s}}}t_j<\tilde{\mathfrak{m}}$ for all $s$, it follows that $I_{\disp{\disp{N_s}}}^d\leq 1-\frac{\disp{r\mathfrak{t}}}{\disp{\tilde
{\mathfrak{m}}}}<1$, for all $s$. Hence $I_N^d\nrightarrow 1$ as $N\rightarrow\infty$, a contradiction.
\end{proof}
 
Let us now concentrate on those orbits whose frequency of visiting the cusp is zero. A complete characterization of such orbits is given by the following proposition.
\begin{proposition}\label{Cusp 3.6}
If $\frac{\displaystyle 1}{\displaystyle N}(\log|a_{\disp{j_1}}|+\log|a_{\disp{j_2}}|.......+
\log|a_{\disp{j_{\disp{\mathfrak{j}_{\underline{d}}^N}}}}|)
\rightarrow 0$ as $N\rightarrow\infty$ for some $d>1$,
then $I_N^{{d}'}=\frac{\disp1}{\disp{S_N}}\disp{\int_{\disp{0}}^{\disp{S_N}}
\chi_{d'}(g_tv)dt}\rightarrow 0$ as $N\rightarrow\infty$
for all $d'>d$. On the other hand if $I_N^{d'}\rightarrow 0$ as $N\rightarrow\infty$ for some 
$d'>1$, then $\frac{\displaystyle 1}{\displaystyle N}(\log|a_{\disp{j_1}}|+\log|a_{\disp{j_2}}|
 .......+\log|a_{\disp{j_{\disp{\mathfrak{j}_{\bar{d}}^N}}}}|)\rightarrow 0$  as 
$N\rightarrow\infty$ for all $d>d'$.
\end{proposition}
\begin{proof}
From Lemma \ref{Cusp 3.3}, we have
\begin{align}\label{equation3}
\frac{\disp1}{\disp N}\disp{\sum\limits_{\disp{k=1}}^{\disp{\mathfrak{j}_ {\underline{d}}^N}}
 2\log\disp{|a_{\disp{j_k}}|}}-\frac{\disp1}{\disp N}\mathfrak{j}_{\underline{d}}^N \kappa
 \leq \frac{\disp1}{\disp N}\disp{\sum\limits_{\disp{k=1}}^{\disp{\mathfrak{j}
 _{\underline{d}}^N}}t_{\disp{\disp{j_k}}}}
 \leq \frac{\disp1}{\disp N}\disp{\sum\limits_{\disp{k=1}}^{\disp{\mathfrak{j}
 _{\underline{d}}^N}}2\log\disp{|a_{\disp{j_k}}|}}+\frac{1}{N}\mathfrak{j}_{\underline{d}}^N \kappa
 \end{align}
 where $\kappa$ is as in that Lemma. Note that  
 $\frac{\disp1}{\disp N}\disp{\sum\limits_{\disp{k=1}}^{\disp{\mathfrak{j}_
 {\underline{d}}^N}}2\log\disp{|a_{\disp{j_k}}|}}\rightarrow 0$ 
 implies $\frac{\disp1}{\disp N}\mathfrak{j}_{\underline{d}}^N\rightarrow 0$ as $N\rightarrow\infty$.
 Then from (\ref{equation3}), we conclude that
 $\frac{\displaystyle 1}{\displaystyle N}(\log|a_{\disp{j_1}}|+\log|a_{\disp{j_2}}|.......+
  \log|a_{\disp{j_{\disp{\mathfrak{j}_{\underline{d}}^N}}}}|)
 \rightarrow 0$  is equivalent to 
 $\frac{\disp1}{\disp N}\disp{\sum\limits_{\disp{k=1}}^{\disp{\mathfrak{j}_
 {\underline{d}}^N}}t_{\disp{\disp{j_k}}}}\rightarrow 0$ as $N\rightarrow\infty$.
 
 Now for any $d'> d$,
 \begin{center}
 $I_N^{d'}\leq I_N^d=\frac{\disp{1}}{\disp{S_N}}\disp{\int_{\disp{0}}^{\disp{S_N}}\chi_d(g_tv)dt}
 \leq \frac{\disp{\frac{1}{N}\sum\limits_{\disp{k=1}}^{\disp{\mathfrak{j}_
 {\underline{d}}^N}}}
\disp{t_{\disp{\disp{j_k}}}}}{\disp{\frac{1}{N}\sum\limits_{\disp{j=1}}^{\disp{N}} t_j}}$,
\end{center}
which tends to $0$ as $N\rightarrow\infty$ since 
$\disp{\frac{1}{N}\sum\limits_{\disp{j=1}}^{\disp N} t_j}$ is bounded below by $\mathfrak{t}$. 

To prove the converse statement, let us assume that $I_N^{{d}'}
\rightarrow 0$ as $N\rightarrow\infty$, and $d>d'$. Suppose 
$\frac{\displaystyle 1}{\displaystyle N}(\log|a_{\disp{j_1}}|+\log|a_{\disp{j_2}}|.......+
 \log|a_{\disp{j_{\disp{\mathfrak{j}_{\bar{d}}^N}}}}|)\nrightarrow 0$ as $N\rightarrow\infty$. Then using another version of (\ref{equation3}), with $\underline{d}$ replaced by $\bar{d}$, there is a subsequence $\{\disp{N_s}\}$ and $r>0$, such that $\disp{\frac{\displaystyle1}{\displaystyle{\disp{N_s}}}\sum\limits_{\disp{k=1}}^
{\disp{\mathfrak{j}_{\bar{d}}^{\disp{N_s}}}}}\disp{t_{\disp{\disp{j_k}}}}>r$ 
for all $s\in\N$. Note that, as $I_{\disp{N_s}}^d\leq I_{\disp{N_s}}^{{d}'}\rightarrow 0$ when $s\rightarrow\infty$, we have $\frac{\disp1}{\disp{\disp{N_s}}}\mathfrak{j}_{\bar{d}}^{\disp{N_s}}\rightarrow 0$ as $s\rightarrow\infty$. Because if $\frac{\disp1}{\disp{\disp{N_s}}}\mathfrak{j}_{
\bar{d}}^{\disp{N_s}}\nrightarrow 0$ as $s\rightarrow\infty$, then $\frac{\disp1}{\disp{\disp{N_s}}}\mathfrak{j}_{\bar{d}}^{\disp{N_s}}>\tilde r$ for some $\tilde r>0$ and for infinitely many $s\in\N$. Then $\disp{\frac{1}{\disp{N_s}}\sum\limits_{\disp{k=1}}^{\disp{\mathfrak{j}_{\bar{d'}}^{\disp{N_s}}}}}
\disp{(t_{\disp{\disp{j_k}}}-c_{\disp{\disp{j_k}}})}>\tilde r c_1$, for infinitely many $s$, which in turn implies that $I_{\disp{N_s}}^{d'}>\tilde{r}c_1$ for infinitely many $s\geq1$, where $t_{\disp{\disp{j_k}}}-c_{\disp{\disp{j_k}}}>c_1>0$. This is a contradiction to the fact that $I_N^{{d}'}\rightarrow 0$ as $N\rightarrow\infty$.  

Now let $h_{a,b}^d$ denote the least upper bound of the distances from the cross-section point on $C$ to the horizontal line $y=d$, for all $(a,b)$-reduced geodesics. Then
\begin{align}\label{equation4}
 I_{\disp{N_s}}^{d'}\geq I_{\disp{N_s}}^d&\geq \frac{\disp{\frac{1}
 {\disp{N_s}}\sum\limits_
 {\disp{k=1}}^{\disp{\mathfrak{j}_{\bar{d}}^{\disp{N_s}}}}}
 \disp{(t_{\disp{\disp{j_k}}}-c_
 {\disp{\disp{j_k}}})}}{\disp{\frac{1}{\disp{N_s}}
 \sum\limits_{\disp{j=1}}^{\disp{N_s}} t_j}}\\
 \label{equation5}&\geq\frac{\disp{\frac{1}{\disp{N_s}}\sum\limits_{\disp{k=1}}^
 {\disp{\mathfrak{j}_
 {\bar{d}}^{\disp{N_s}}}}}\disp{t_{\disp{\disp{j_k}}}-\frac{1}{\disp{N_s}}
 \disp{\mathfrak{j}_{\bar{d}}^{\disp{N_s}}
 2h_{a,b}^d}}}{\disp{\frac{1}{\disp{N_s}}\sum\limits_{\disp{j=1}}^{\disp{N_s}} t_j}}.
\end{align}
Since $\frac{\disp1}{\disp{\disp{N_s}}}\mathfrak{j}_{\bar{d}}
^{\disp{N_s}}\rightarrow 0$ as $s\rightarrow\infty$, it follows that  there is some $\mathfrak{m}>0$ such that\\ $\disp{\frac{1}{\disp{N_s}}\sum\limits_{\disp{j=1}}^{\disp{N_s}} t_j}
<\mathfrak{m}$ for all $s\in\N$. 
Therefore, from (\ref{equation4}) and (\ref{equation5}), we conclude that there exists some $0<r_1<r$, such that $I_{\disp{N_s}}^{d'}>\frac{\disp{r_1}}
{\disp{\mathfrak{m}}}$ for sufficiently large $s$. Which is a contradiction to the assumption that $I_{N}^{d'}
\rightarrow 0$ as $N\rightarrow\infty$. This completes the proof of the proposition. 
\end{proof}

Now the proof of Theorem \ref{Cusp 1.1} for $(a,b)\in\mathcal{S}'$ follows from
Proposition \ref{Cusp 3.5} and Proposition \ref{Cusp 3.6}.
\begin{remark}
Note that in Proposition \ref{Cusp 3.5} and Proposition \ref{Cusp 3.6}, we have considered an
$(a,b)$-reduced lift of $\gamma_v$, whereas in Theorem \ref{Cusp 1.1}, we have considered 
any lift of $\gamma_v$ to $\H$. This does not lead to any ambiguity because  
if we obtain an $(a,b)$-reduced geodesic $\gamma$ with attracting end point 
$w=[a_0,a_1,a_2,...]_{a,b}$, from a geodesic $\gamma'$ with attracting end point $w'=[a'_0,a'_1,
a'_2,...]_{a,b}$, then $a_j=a'_{j+n}$ for some $n\in\N$. 
\end{remark}

\subsection{$(a,b)=(-1,1)$}
Now we concentrate on the special case $(a,b)=(-1,1)$. Recall that the coding of geodesics on the modular surface using this particular continued fraction is discussed in detail in \cite{KU}, where it is called the alternating continued fraction coding. The name alternating continued fraction comes from the fact that the partial quotients of the $(-1,1)$-continued fraction expansion
of a real number has alternate signs. This particular coding procedure does not provide a cross-section contained in a compact subset of $T^1M$. Recall form \cite{KU}, that a geodesic in $\H$ is called $A$-reduced ($(-1,1)$-reduced with our convention), if its attracting endpoint $w$ and repelling endpoint $u$ satisfy $|w|>1$ and $-1<\text{sgn}(w)u<0$ respectively. So the cross-section point for an $A$-reduced geodesic can be as close to the real line as one wants, showing that the cross-section is not contained inside a compact set in $T^1M$. So the $j$th return time may not be at a bounded distance from $2\log|a_j|$. But $t_j$ can be controlled using a couple of preceding and couple of succeeding entries in the sequence of partial quotients. We recall from \cite{KU}, the following formula for the $j$th return time: 
\begin{align*}
 t_j=2\log |w_j|+\log\frac{\displaystyle{|w_j-u_j|\sqrt{w_j^2-1}}}{\displaystyle{w_j^2
 \sqrt{1-u_j^2}}}-\log\frac{\displaystyle{|w_{j+1}-u_{j+1}|\sqrt{w_{j+1}^2-1}}}
 {\displaystyle{w_{j+1}^2\sqrt{1-u_{j+1}^2}}}.
\end{align*}
Now assume that $w_j>0$, then it follows from the definition of $A$-reduced geodesics that $u_j<0$. Since the partial quotients have alternate signs, we also have $w_{j+1}<0$ and consequently $u_{j+1}> 0$. Then,
\begin{align*}
 t_j\leq 2\log |w_j|+\left|\log\frac{|\displaystyle{1-\frac{u_j}{w_j}}|}{|\displaystyle{1-
 \frac{u_{j+1}}{w_{j+1}}}|}\right|+\\ \frac{\displaystyle{1}}{\displaystyle{2}}\left[
 \left|\log\left(1-\frac{1}{w_j}\right)\right|+|\log(1+u_j)|+\left|\log\left(1+\frac{1}{w_{j+1}}
 \right)\right|+
 |\log(1-u_{j+1})|\right]+\\ \frac{\displaystyle{1}}{\displaystyle{2}}\left|\log
 \frac{\displaystyle{\left(1+\frac{1}{w_j}\right)(1+u_{j+1})}}{\displaystyle{(1-u_j)
 \left(1-\frac{1}{w_{j+1}}\right)}}\right|.
\end{align*}
\begin{figure}\label{hd1}
\centering
\thispagestyle{empty}
\includegraphics[scale=0.7]{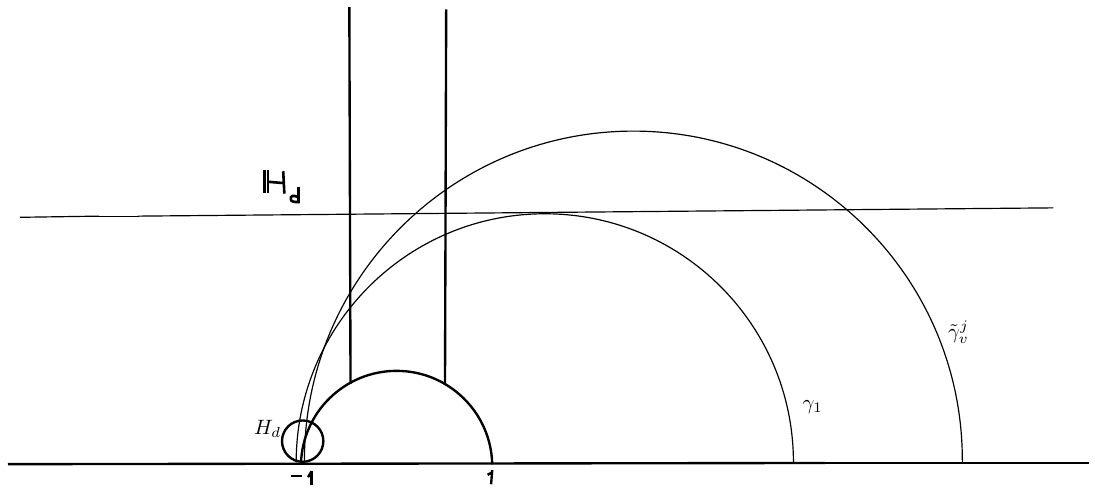}
\caption{Possible connected components of $\pi^{-1}(\gamma_v^j\cap M_d)$}
\end{figure}
Now using the assumption that $a_j>0$ and consequently $a_{j+1}<0$, $a_{j+2}>0$, it is easy to see that $$1-\frac{\displaystyle{1}}{\displaystyle{w_{j}}}\geq 
1-\frac{\displaystyle{1}}{\displaystyle{1+\frac{\displaystyle{1}}{\displaystyle{|a_{j+1}|+\delta}}}},$$ where $\delta$ is some real number such that $0\leq\delta\leq1$. Then it follows that, $$\left|\log\left(1-\frac{\displaystyle{1}}{\displaystyle{w_j}}\right)\right|
\leq\log|a_{j+1}|+\log3.$$ By a similar reasoning, $$1+u_j\geq 1-\frac{\displaystyle{1}}{\displaystyle{1+\frac{\displaystyle{1}}{\displaystyle{|a_{j-2}|+\delta^{\prime}}}}},$$ with $0\leq\delta^{\prime}\leq1$, and it follows that, $$|\log(1+u_j)|\leq\log|a_{j-2}|+\log3.$$
Using the continued fraction expansions for $w_{j+1}$ and $u_{j+1}$, we obtain similar estimates for other quantities in the above inequality involving $t_j$. The case $w_j<0$ can be treated similarly and we get the following estimate for the return time $t_j$:
\begin{align}\label{equation6}
 t_j\leq 2\log|a_j|+2\hspace*{0.1cm}\text{max}\{\log|a_{j+1}|+\log|a_{j-1}|+\log|a_{j+2}|+
 \log|a_{j-2}|\}+c, 
\end{align}
here $c$ is some constant which is independent of $j$. On the other hand, considering the definition of $A$-reduced geodesics, and the fact that the length of the geodesic segment joining the point $i$ 
and $k+i$ is at a bounded distance from $2\log |k|$, independent of $k\in \Z$, it is easy to see that 
\begin{align}\label{equation7}
 t_j\geq 2\log|a_j|-c',
\end{align}
where $c'$ can be taken as the hyperbolic length of the segment of the unit circle joining the point $i$ and $\frac{\displaystyle{1}}{\displaystyle{2}}+\frac{\displaystyle{\sqrt{3}}}
{\displaystyle{2}}i$.

Also note that in this special case, whenever $\gamma_v^j\cap M_d$ is non-empty, the number of connected components of $\pi^{-1}(\gamma_v^j\cap M_d)$ can be more than one, in fact it can be at most three. One component is $\tilde{\gamma}_v^j\cap\H_d$; one of the other two may be the segment starting from the cross-section point up to the 
intersection point of $\tilde{\gamma}_v^j$ with the horocycle $H_d$, where $H_d$ is the image of the horocycle $y=d$ under $T^{-1}S$ as shown in Figure $4$; the third component may be a similar one coming from near the other end of $\tilde{\gamma}_v^j$. In Figure $4$, the geodesic $\gamma_1$ is the geodesic which is tangent to the horizontal line $y=d$ and passes through the intersection point of the vertical line based at $-1$ and the horocycle $H_d$. Let $h_u^d$ be the hyperbolic length of the segment of $\gamma_1$ joining the pair of points where it cuts the horocycle $H_d$ and where it touches the line $y=d$. Then $h(\gamma_v^j
\backslash M_d)\leq 2 h_u^d$. On the other hand, let $h_l^d$ denote the hyperbolic distance between the points $i$ and the horizontal line $y=d$. Then if $\gamma_v^j\cap M_d\neq\phi$, $h(\gamma_v^j\backslash M_d)\geq 2 h_l^d$. Using these observations, the following two propositions from which the proof of Theorem \ref{Cusp 1.1} follows in the case $(a,b)=(-1,1)$, can be proved by adopting the similar strategies as in the proofs of Proposition \ref{Cusp 3.5} and Proposition \ref{Cusp 3.6} respectively.

\begin{proposition}\label{Cusp 4.1}
Let $v\in T^1M$, $\gamma_v$ be the corresponding geodesic on $M$ and $\tilde{\gamma}_v$ be an $A$-reduced lift of $\gamma_v$ in $\H$. Let $w=[a_0,a_1,a_2,...]_{(-1,1)}$ be the attracting endpoint of $\tilde{\gamma}_v$. Then
$I_N^d=\frac{\disp{1}}{\disp{S_N}}\disp{\int_{\disp{0}}^{\disp{S_N}}\chi_d(g_tv)dt}\rightarrow 1$ as $N\rightarrow\infty$ for all $d>1$ if and only if \\ ${\frac{\displaystyle1}
{\displaystyle{N}}\displaystyle{\sum\limits_{\disp{j=1}}^{\disp{N}} t_j}}\rightarrow\infty$ as $N\rightarrow\infty$.
\end{proposition}

\begin{proposition}\label{Cusp 4.2}
If $\frac{\displaystyle1}{\displaystyle N}\displaystyle{\sum\limits_{\disp{t_j>\mathfrak{c}, 1\leq j\leq N}}t_j}\rightarrow 0$ as $N\rightarrow\infty$ for some $\mathfrak{c}>0$, 
then there exists $d>1$, such that $I_N^{d'}\rightarrow 0$ as $N\rightarrow
\infty$ for all $d'> d$. On the other hand, if 
$I_N^d\rightarrow 0$ as $N\rightarrow\infty$ for some $d>1$, then there exist $\mathfrak{c}>0$
such that $\frac{\displaystyle1}{\displaystyle N}\displaystyle{\sum\limits_{\disp{t_j>
\mathfrak{c}',1\leq j\leq N}}
t_j}\rightarrow 0$ as $N\rightarrow\infty$ for all $\mathfrak{c}'> \mathfrak{c}$.
\end{proposition}

{\bf Proof of Theorem \ref{Cusp 1.1} in the case of $(-1,1)$-continued fraction.}\\
It follows easily from (\ref{equation6}) and (\ref{equation7}), that $\frac{\displaystyle{1}}
{\displaystyle{N}}\disp{\sum\limits_{\disp{j=0}}^{\disp{N-1}}\log|a_j|}\rightarrow\infty$ as
$N\rightarrow\infty$ is equivalent to $\frac{\displaystyle{1}}{\displaystyle{N}}
\disp{\sum\limits_{\disp{j=1}}^{\disp{N}} t_j}\rightarrow\infty$ as $N\rightarrow\infty$. 
Now let $$\mathfrak{j}_d^N=\#\{1\leq j\leq N:\gamma_v^j\cap M_d\neq\phi\}.$$
Then for $d>1$, from (\ref{equation6}) we get 
\begin{align}\label{equation8}
\displaystyle{\sum\limits_{\displaystyle{t_j>10\log d,1\leq j\leq N}}t_j}
\leq\displaystyle{\sum\limits_{\displaystyle{|a_j|>d,-2\leq j
\leq(N+2)}}10\log|a_j|}+\mathfrak{j}_d^N c.
\end{align}
Therefore, if $\frac{\displaystyle{1}}{\displaystyle{N}}\displaystyle{\sum
\limits_{\displaystyle{|a_j|>d,0\leq j\leq N-1}}\log|a_j|}\rightarrow 0$ as $N\rightarrow\infty$ for some $d>1$, which also implies $\frac{\displaystyle{1}}{\displaystyle{N}}
\mathfrak{j}_d^N\rightarrow 0$ as $N\rightarrow\infty$, then it follows from (\ref{equation8}), that\\  $\frac{\displaystyle{1}}{\displaystyle{N}}\displaystyle{\sum
\limits_{\displaystyle{t_j>d',1\leq j\leq N}}t_j}\rightarrow 0$ as $N\rightarrow\infty$ 
for all $d'>10\log d$.
On the other hand, it follows easily from (\ref{equation7}), that if 
$\frac{\displaystyle{1}}{\displaystyle{N}}\displaystyle{\sum
\limits_{\displaystyle{t_j>d,1\leq j\leq N}}t_j}\rightarrow 0$ for some $d>0$, which also 
implies $\frac{\displaystyle{1}}{\displaystyle{N}}
\mathfrak{j}_d^N\rightarrow 0$ as $N\rightarrow\infty$, there 
exists $d'>1$, such that\\ $\frac{\displaystyle{1}}{\displaystyle{N}}\displaystyle{\sum
\limits_{\displaystyle{|a_j|>d'',0\leq j\leq N-1}}\log|a_j|}\rightarrow 0$ as $N\rightarrow
\infty$ for all $d''> e^{\disp{d'/2}}$. With these observations, now the proof of Theorem \ref{Cusp 1.1} 
in the case of $(-1,1)$-continued fraction, follows from Proposition \ref{Cusp 4.1} and 
Proposition \ref{Cusp 4.2}. \hspace*{5.5cm}  $\Box$

\section{Acknowledgements} The author is thankful to S. G. Dani for suggesting the problem and his constant help in writing the paper. Thanks are also due to the referee of this version for valuable suggestions which have helped to improve the exposition of the article. The author thanks Indian Statistical Institute Bangalore, Harish Chandra Research Institute Allahabad and Indian Statistical Institute Kolkata for their hospitality during the author's stay there, which has made this work possible. Financial support from the National Board for Higher Mathematics India through NBHM post-doctoral fellowship is duly acknowledged.

\smallskip
\bibliography{ref}
\bibliographystyle{plain}

\end{document}